\documentclass[reqno]{amsart}
\usepackage{amssymb}
\usepackage[pagebackref,hypertex]{hyperref}
\date{This manuscript was completed on 10 January 2008}
\date{}
\allowdisplaybreaks[4]
\theoremstyle{plain}
\newtheorem{thm}{Theorem}
\theoremstyle{remark}
\newtheorem{rem}{Remark}
\newcommand{\td}{\textup{d}}

\begin{document}

\title[Proofs for a monotonicity result of a function]
{Proofs for a monotonicity result of a function involving the psi and exponential functions}

\author[F. Qi]{Feng Qi}
\address[F. Qi]{Department of Mathematics, College of Science, Tianjin Polytechnic University, Tianjin City, 300160, China}
\email{\href{mailto: F. Qi <qifeng618@gmail.com>}{qifeng618@gmail.com}, \href{mailto: F. Qi <qifeng618@hotmail.com>}{qifeng618@hotmail.com}, \href{mailto: F. Qi <qifeng618@qq.com>}{qifeng618@qq.com}}
\urladdr{\url{http://qifeng618.wordpress.com}}

\author[B.-N. Guo]{Bai-Ni Guo}
\address[B.-N. Guo]{School of Mathematics and Informatics,
Henan Polytechnic University, Jiaozuo City, Henan Province, 454010, China} \email{\href{mailto: B.-N.
Guo <bai.ni.guo@gmail.com>}{bai.ni.guo@gmail.com}, \href{mailto: B.-N. Guo
<bai.ni.guo@hotmail.com>}{bai.ni.guo@hotmail.com}}

\begin{abstract}
In the note, two alternative proofs are provided for a monotonicity result that the function $\psi(x)+\ln\bigl(e^{1/x}-1\bigr)$ is strictly increasing on $(0,\infty)$, where $\psi(x)$ is the psi function.
\end{abstract}


\thanks{The first author was supported in part by the China Scholarship Council and the Science Foundation of Tianjin Polytechnic University}

\maketitle

\section{Introduction}

It is well-known that the classical gamma function
\begin{equation}\label{gamdef}
\Gamma(x)=\int^\infty_0t^{x-1} e^{-t}\td t
\end{equation}
for $x>0$, the psi function $\psi(x)=\frac{\Gamma'(x)}{\Gamma(x)}$, and the polygamma functions $\psi ^{(k)}(x)$ for $k\in\mathbb{N}$ play central roles in the theory of special functions and have much extensive applications in many branches.
\par
In~\cite[Theorem~2.1]{batir-new}, it was discovered that if $a\le-\ln2$ and $b\ge0$, then
\begin{equation}\label{batir-ineq-orig}
a-\ln\bigl(e^{1/x}-1\bigr)<\psi(x)<b-\ln\bigl(e^{1/x}-1\bigr)
\end{equation}
holds for $x>0$.
\par
In~\cite[pp.~386--388]{alzer-expo-math-2006}, the function
\begin{equation}
\phi(x)=\psi(x)+\ln\bigl(e^{1/x}-1\bigr)
\end{equation}
was proved to be strictly increasing on $(0,\infty)$ and
\begin{equation}\label{limit=0}
\lim_{x\to\infty}\phi(x)=0.
\end{equation}
\par
In~\cite[Theorem~2.8]{batir-jmaa-06-05-065}, the inequality~\eqref{batir-ineq-orig} was sharpened as follows: If and only if $a\le-\gamma$ and $b\ge0$, the inequality~\eqref{batir-ineq-orig} is valid for $x>0$, where $\gamma=0.577\dotsc$ stands for Euler-Mascheroni's constant.
\par
In~\cite{property-psi.tex} and its preprint~\cite{property-psi.tex-arXiv}, among other things, the function $\phi(x)$ was proved to be not only strictly increasing but also strictly concave on $(0,\infty)$, with~\eqref{limit=0} and
$$
\lim_{x\to0^+}\phi(x)=-\gamma.
$$
\par
In~\cite{alzer-expo-math-2006, property-psi.tex, property-psi.tex-arXiv}, proofs of the monotonicity result that the function $\phi(x)$ is strictly increasing occupy almost two printed pages respectively.
\par
The aim of this note is to provide two alternative proofs for the monotonicity of the function $\phi(x)$, which may be recited as Theorem~\ref{peoperty-psi-ii.tex-thm} below.

\begin{thm}\label{peoperty-psi-ii.tex-thm}
The function $\phi(x)$ is strictly increasing on $(0,\infty)$.
\end{thm}

\begin{rem}
The above Theorem~\ref{peoperty-psi-ii.tex-thm} complements the well-known Bohr-Mollerup theorem which says that the gamma function is logarithmically convex, that is, the psi function is increasing.
\end{rem}

\begin{rem}
It is well-known that the $n$-th harmonic number is defined by
\begin{equation}
H_n=\sum_{k=1}^n\frac1k
\end{equation}
for $n\in\mathbb{N}$ and that $H_n$ can be expressed in terms of the psi function $\psi(x)$ by
\begin{equation}
H_n=\psi(n+1)+\gamma.
\end{equation}
Consequently, the monotonicity result of $\phi(x)$ implies the sharp double inequality in~\cite[Theorem~2.8]{batir-jmaa-06-05-065} and the sharp inequalities for harmonic numbers $H_n$ in~\cite[pp.~386\nobreakdash--387]{alzer-expo-math-2006}: For $n\in\mathbb{N}$, we have
\begin{equation}\label{aler-harmonic-ineq}
1+\ln\bigl(\sqrt{e}\,-1\bigr)\le H_n+\ln\bigl(e^{1/(n+1)}-1\bigr)<\gamma.
\end{equation}
The constants $1+\ln\bigl(\sqrt{e}\,-1\bigr)$ and $\gamma$ in~\eqref{aler-harmonic-ineq} are the best possible.
\par
Some sharp inequalities for the $n$-th harmonic numbers were also established in~\cite{Chao-Ping-AML-10, ucdavis-harmonic, harseq, harmonic-number-refine-chen.tex-conf, qi-cui-jmaa, harmonic-number-refine-chen.tex, Infinite-family-Digamma.tex} and closely-related references therein. Especially, we would like to mention that a new class of sequences of the form
\begin{equation}
\mu_n=H_n+\ln\bigl(e^{a/(n+b)}-1\bigr)
\end{equation}
was introduced in~\cite{Mortici-10-Carpathian} and the fastest sequence $\{\mu_n\}_{n\ge1}$ was obtained for $a=\frac1{\sqrt2\,}$ and $b=\frac{2+\sqrt2\,}{4}$.
\end{rem}

\section{Two alternative proofs of Theorem~\ref{peoperty-psi-ii.tex-thm}}

Now we are in a position to provide two alternative proofs of Theorem~\ref{peoperty-psi-ii.tex-thm}.

\begin{proof}[First proof]
It is well-known that
\begin{equation}\label{recursion-gamma}
\Gamma(x+1)=x\Gamma(x)
\end{equation}
for $x>0$. Taking the logarithm and differentiating on both sides of~\eqref{recursion-gamma} yield
\begin{equation}\label{psisymp4}
\psi^{(i-1)}(x+1)-\psi^{(i-1)}(x)=(-1)^{i-1}\frac{(i-1)!}{x^i}, \quad i\in\mathbb{N}.
\end{equation}
Therefore, the function
\begin{equation}\label{psi(x)2+psi(x)-fun}
h(x)=[\psi'(x)]^2+\psi''(x)
\end{equation}
on $(0,\infty)$ satisfies
\begin{align*}
  h(x)-h(x+1)&=[\psi'(x)-\psi'(x+1)][\psi'(x)+\psi'(x+1)]+\psi''(x)-\psi''(x+1)\\
  &=\frac2{x^2}\biggl[\psi'(x)-\frac1x-\frac1{2x^2}\biggr]\\
  &\triangleq\frac2{x^2}g(x)
\end{align*}
and the difference
\begin{equation*}
g(x+1)-g(x)=\psi'(x+1)-\psi'(x)+\frac{2 x^2+4 x+1}{2 x^2 (x+1)^2}=-\frac{1}{2 x^2 (x+1)^2}<0.
\end{equation*}
Hence, by induction, we have
\begin{align*}
g(x)>g(x+1)>g(x+2)>\dotsm>g(x+k)\to0
\end{align*}
as $k\to\infty$, where we utilize the well-known fact that $\psi^{(i)}(x)\to0$ as $x\to\infty$ for $i\in\mathbb{N}$. By the same argument as the deduction of $g(x)>0$, we can obtain $h(x)>0$ on $(0,\infty)$. Consequently, the exponential function of $\phi(x)$ satisfies
\begin{gather*}
e^{\phi(x)}=e^{\psi(x)}\bigl(e^{1/x}-1\bigr)=e^{\psi(x)+1/x}-e^{\psi(x)} =e^{\psi(x+1)}-e^{\psi(x)}\triangleq f(x),\\
f'(x)=e^{\psi(x+1)}\psi'(x+1)-e^{\psi(x)}\psi'(x) \triangleq q(x+1)-q(x),\\
q'(x)=[e^{\psi(x)}\psi'(x)]'=e^{\psi(x)}\bigl\{\psi''(x)+[\psi'(x)]^2\bigr\}=e^{\psi(x)}h(x)>0.
\end{gather*}
As a result, the function $q(x)$ is strictly increasing, and so $f'(x)>0$ on $(0,\infty)$. As a result, the function $f(x)$, and then $\phi(x)$, is strictly increasing on $(0,\infty)$. The first proof of Theorem~\ref{peoperty-psi-ii.tex-thm} is complete.
\end{proof}

\begin{rem}
In~\cite[p.~208]{forum-alzer} and~\cite[Lemma~1.1]{batir-new}, the positivity of the function $h(x)$ defined by \eqref{psi(x)2+psi(x)-fun} on $(0,\infty)$ was verified by different approaches. Recently this positivity was further generalized to a complete monotonicity result in \cite{Guo-Qi-Srivasta-Unique.tex, notes-best-simple-open-jkms.tex, x-4-di-tri-gamma-p(x).tex, notes-best-simple-cpaa.tex, AAM-Qi-09-PolyGamma.tex} and closely-related references therein.
\end{rem}

\begin{rem}
The first proof of Theorem~\ref{peoperty-psi-ii.tex-thm} shows that the function $e^{\psi(x)}$ is convex on $(0,\infty)$. It is common knowledge that the psi function $\psi(x)$ is concave on $(0,\infty)$. This gives an example that a logarithmically concave function is not concave but convex. For more information, please see~\cite[Section~3]{merkle-via} and~\cite[p.~6, Remark~1.9]{bounds-two-gammas.tex}.
\end{rem}

\begin{proof}[Second proof]
This is based on the observation that $f'(x)>0$ is equivalent to
\begin{equation}\label{Mortici-sugg-1}
  e^{1/x}\psi'(x+1)>\psi'(x)
\end{equation}
for $x>0$.
\par
The inequality
\begin{equation}\label{Mortici-sugg-2}
  \frac1x+\frac1{2x^2}+\frac1{6x^3}-\frac1{30x^5}<\psi'(x)<\frac1x+\frac1{2x^2}+\frac1{6x^3}
\end{equation}
for $x>0$ arises from the standard asymptotic series of $\psi'(x)$, which was proved in~\cite{Mortici-10-Asymptot-Anal} to be of great help in establishing other results.
\par
Taking into account \eqref{Mortici-sugg-2}, in order to prove~\eqref{Mortici-sugg-1}, it suffices to show
\begin{equation*}
e^{1/x}\biggl[\frac1{x+1}+\frac1{2(x+1)^2}+\frac1{6(x+1)^3}-\frac1{30(x+1)^5}\biggr] >\frac1x+\frac1{2x^2}+\frac1{6x^3}
\end{equation*}
or, equivalently,
\begin{align*}
p(x)&=\frac1x+\ln\biggl[\frac1{x+1}+\frac1{2(x+1)^2}+\frac1{6(x+1)^3}-\frac1{30(x+1)^5}\biggr]\\
&\quad -\ln\biggl(\frac1x+\frac1{2x^2}+\frac1{6x^3}\biggr)\\
&>0.
\end{align*}
But
\begin{equation*}
p'(x)=-\frac{49+224x+475x^2+521x^3+246x^4+45x^5}{x^2(x+1)(6x^2+3x+1)(49+175x+230x^2+135x^3+30x^4)}<0.
\end{equation*}
Finally, the function $q(x)$ is strictly decreasing with $\lim_{x\to\infty}q(x)=0$, so $q(x)>0$. The second proof is complete.
\end{proof}

\begin{rem}
We note that the second proof of Theorem~\ref{peoperty-psi-ii.tex-thm} is suggested to be added by an anonymous referee of this paper.
\end{rem}

\begin{rem}
This is a revised and expanded version of the preprint~\cite{property-psi-ii.tex}.
\end{rem}

\subsection*{Acknowledgements}
The authors heartily appreciate the anonymous referees for their crucial comments and valuably technical advices on the original version of this paper.

\end{document}